\documentclass[12pt,leqno]{article}
\usepackage{amsfonts,amsthm,amsmath}
\newtheorem{thm}{Theorem}
\newtheorem{prop}[thm]{Proposition}

\theoremstyle{remark}

\newcommand{\ZZ}{\mathbb{Z}}
\newcommand{\RR}{\mathbb{R}}

\DeclareMathOperator{\Aut}{Aut}

\DeclareMathOperator{\sub}{sub}

\begin{document}

\title{Optimal Self-Dual $\ZZ_4$-Codes and 
a Unimodular Lattice in Dimension $41$}

\author{Masaaki Harada\thanks{Department of Mathematical Sciences,
Yamagata University, Yamagata 990--8560, Japan, and
PRESTO, Japan Science and Technology Agency, Kawaguchi,
Saitama 332--0012, Japan. email: mharada@sci.kj.yamagata-u.ac.jp}}

\maketitle

\begin{abstract}
For lengths up to $47$ except $37$,
we determine the largest minimum Euclidean weight 
among all Type~I $\ZZ_4$-codes of that length.
We also give the first example of an optimal
odd unimodular lattice in dimension $41$
explicitly, which is constructed from some 
Type~I $\ZZ_4$-code of length $41$.
\end{abstract}

\section{Introduction}\label{Sec:1}

Let $\ZZ_4\ (=\{0,1,2,3\})$ denote the ring of integers
modulo $4$.
A $\ZZ_{4}$-code $C$ of length $n$ 
is a $\ZZ_{4}$-submodule of $\ZZ_{4}^n$.
The dual code $C^\perp$ of $C$ is defined as
$\{ x \in \ZZ_{4}^n \mid x \cdot y = 0$ for all $y \in C\}$
under the standard inner product~$x \cdot y$.
A code $C$ is {\em self-dual} if $C=C^\perp.$
The Euclidean weight of a codeword $x=(x_1,\ldots,x_n)$ is
$n_1(x)+4n_2(x)+n_3(x)$, where $n_{\alpha}(x)$ denotes 
the number of components $i$ with $x_i=\alpha$ $(\alpha=1,2,3)$.
The minimum Euclidean weight $d_E(C)$ of $C$ is the smallest Euclidean
weight among all nonzero codewords of $C$.
A self-dual code which has the property that all
Euclidean weights are divisible by eight,
is called {\em Type~II} \cite{Z4-BSBM}
(see also \cite{Z4-HSG}).
A self-dual code which is not Type~II, is 
called {\em Type~I}.
A Type~II $\ZZ_4$-code of length $n$ exists if and 
only if $n \equiv 0 \pmod 8$~\cite{Z4-BSBM},
while a Type~I $\ZZ_4$-code exists for every length.

It was shown in~\cite{Z4-BSBM} that the minimum Euclidean 
weight $d_E(C)$ of a Type~II code $C$ 
of length $n$ is bounded by
$d_E(C) \le 8 \lfloor \frac{n}{24} \rfloor +8$.
A Type~II code meeting this bound with equality is called 
{\em extremal}.
It was also shown in \cite{Rains-Sloane98} that
the minimum Euclidean 
weight $d_E(C)$ of a Type~I code $C$ 
of length $n$ is bounded by
$d_E(C) \le 8 \lfloor \frac{n}{24} \rfloor +8$
if $n \not\equiv 23 \pmod{24}$, and
$d_E(C) \le 8 \lfloor \frac{n}{24} \rfloor +12$
if $n \equiv 23 \pmod{24}$.
It is a fundamental problem to determine the 
largest minimum Euclidean weight among self-dual
codes of that length.
We denote 
the largest minimum Euclidean weight among 
Type~I codes of length $n$ by $d_{max,E}^I(n)$.
These values $d_{max,E}^I(n)$ have been determined
in  \cite{DHS} and \cite{Rains} for $n \le 24$. 
We say that a Type~I code of length $n$ is 
{\em optimal} or {\em Euclidean-optimal}
if it has minimum Euclidean
weight $d_{max,E}^I(n)$.
We pay attention to the Euclidean weight 
from the viewpoint of a connection with
optimal odd unimodular lattices.

In this paper, we determine
the largest minimum Euclidean weight 
$d_{max,E}^I(n)$ for lengths $n \le 47,\ n\not=37$.
To do this, we slightly improve upper bounds on the 
minimum Euclidean weights for lengths 
$n=25$, $26,\ldots,31$, $33,34,35$ (see the bound (\ref{Eq:dE})),
and we construct Type~I $\ZZ_4$-codes
meeting the bound (\ref{Eq:dE}) with equality.
The values $d_{max,E}^I(n)$ are listed in Table \ref{Tab:dE}.
For length $37$, our extensive search failed to discover 
a Type~I $\ZZ_4$-code with minimum Euclidean weight $16$.
However, we have found a Type~I code 
with minimum Euclidean weight $12$.
We also give the first explicit example of an optimal
odd unimodular lattice in dimension $41$, which is constructed 
from an optimal code of length $41$
by Construction A.
All computer calculations in this paper
were done using {\sc Magma} \cite{Magma}.

\begin{table}[thb]
\caption{Largest minimum Euclidean weights of Type~I $\ZZ_4$-codes}
\label{Tab:dE}
\begin{center}
{\small
\begin{tabular}{c|c|c||c|c|c}
\noalign{\hrule height0.8pt}
Length $n$ & $d_{max,E}^I(n)$ & Reference &
Length $n$ & $d_{max,E}^I(n)$ & Reference \\
\hline

25 &  8     & $\sub(C_{26})$  & 37 & 12, 16 &          \\
26 & 12     & $C_{26}$    & 38 & 16     & \cite{GNS}         \\
27 & 12     & $C_{27}$    & 39 & 16     & \cite{GH-Z4-39} \\
28 & 12     & $C_{28}$    & 40 & 16     & \cite{H-Z4-40}  \\
29 & 12     & $C_{29}$    & 41 & 16     & $C_{41}$ \\
30 & 12     & \cite{GNS}  & 42 & 16     & \cite{GNS} \\
31 & 12     & $\sub(C_{32})$ & 43 & 16     & $C_{43}$ \\
32 & 16     & $C_{32}$    & 44 & 16     & $C_{44}$ \\
33 & 12     & $C_{33}$    & 45 & 16     & $C_{45}$ \\
34 & 12     & $C_{34}$    & 46 & 16     & \cite{H-Z4-46}  \\
35 & 12     & $\sub(C_{36})$ & 47 & 16     & \cite{H-Z4-46}  \\
36 & 16     & $C_{36}$       & &&\\
\noalign{\hrule height0.8pt}
   \end{tabular}
}
\end{center}
\end{table}

\section{Preliminaries}

\subsection{Self-dual $\ZZ_4$-codes}

Every $\ZZ_4$-code $C$ of length $n$ has two binary codes 
$C^{(1)}$ and $C^{(2)}$ associated with $C$:
\[
C^{(1)}= \{ c \bmod 2 \mid  c \in C \} \text{  and }
C^{(2)}= \left\{ c \bmod 2 \mid c \in \ZZ_4^n, 2c\in C \right\}.
\]
The binary codes $C^{(1)}$ and $C^{(2)}$ are called the 
{\em residue} and {\em torsion} codes of $C$, respectively.
If $C$ is a self-dual $\ZZ_4$-code, then $ C^{(1)}$ is a binary 
doubly even code with 
$C^{(2)} = {C^{(1)}}^{\perp}$~\cite{Z4-CS}.
It is easy to see that 
\begin{equation}\label{eq:dE}
\min\{d(C^{(1)}),4d(C^{(2)}))\} \le d_E(C)  \le 4d(C^{(2)}),
\end{equation}
where $d(C^{(i)})$ denotes the minimum weight of $C^{(i)}$
$(i=1,2)$.

Codes differing by only a permutation of coordinates are called
permutation-equivalent.
Any self-dual $\ZZ_4$-code of length $n$
is permutation-equivalent to a code 
$C$ with generator matrix 
of the standard form
\begin{equation}
\label{Z4-g-matrix}
\left(\begin{array}{ccc}
I_{k_1} & A & B_1+2B_2 \\
O    &2I_{k_2} & 2D \\
\end{array}\right),
\end{equation}
where $A$, $B_1$, $B_2$ and $D$ are $(1,0)$-matrices,
$I_k$ denotes the identity matrix of order $k$,
and $O$ denotes the zero matrix \cite{Z4-CS}.
The residue code $C^{(1)}$ of $C$ is an $[n,k_1]$ code with 
generator matrix
$
\left(\begin{array}{ccc}
I_{k_1} & A & B_{1}
      \end{array}\right),
$
and the torsion code $C^{(2)}$ is an $[n,k_1+k_2]$ code 
with generator matrix
$
\left(\begin{array}{ccc}
I_{k_1} & A & B_{1}\\
O    &I_{k_2} & D \\
\end{array}\right)
$.

\subsection{Unimodular lattices and upper bounds}
A (Euclidean) lattice $L$ in dimension $n$ is {\em unimodular} if
$L = L^{*}$, where
the dual lattice $L^{*}$ of $L$ is defined as 
$\{ x \in {\RR}^n \mid (x,y) \in \ZZ \text{ for all }
y \in L\}$ under the standard inner product $(x,y)$.
The norm of a vector $x$ is $(x, x)$.
Two lattices $L$ and $L'$ are {\em isomorphic}, denoted $L \cong L'$,
if there exists an orthogonal matrix $A$ with
$L' = L \cdot A =\{xA \mid x \in L\}$.
The minimum norm of $L$ is the smallest 
norm among all nonzero vectors of $L$.
The theta series $\theta_{L}(q)$ of $L$ is the formal power
series
$
\theta_{L}(q) = \sum_{x \in L} q^{(x,x)}
= \sum_{m=0}^{\infty} N_m q^m,
$
where $N_m$ is the number of vectors of norm $m$.
The kissing number is the second nonzero coefficient of the
theta series, that is, the number of vectors of $L$ with
minimum norm.

Let $\mu^O_{max}(n)$  denote the largest minimum
norm among odd unimodular lattices in dimension $n$.
We say that an odd unimodular lattice is
{\em optimal} if it has the largest minimum norm $\mu^O_{max}(n)$.
These values $\mu^O_{max}(n)$ have been determined
for $n \le 47, n \not=37,41$
(see \cite{CS98} and \cite{S-http}).
For $25 \le n \le 47$, the following values are known
\begin{equation}\label{Eq:mu}
\mu^O_{max}(n)=
\left\{\begin{array}{ll}
2    &(n=25),  \\
3    &(n=26,27,\ldots,31,33,34,35), \\
4    &(n=32,36,38,39,40,42,\ldots,47), \\
3 \text{ or } 4 &(n=37,41).
\end{array}\right.
\end{equation}
In this paper, we give the first example of an
odd unimodular lattice in dimension $41$ having
minimum norm $4$.
Hence, we have $\mu^O_{max}(41)=4$.

Let $C$ be a Type~II (resp.\ Type~I) 
$\ZZ_4$-code of length $n$
and minimum Euclidean weight $d_E$.
Then the following lattice
\[
A_{4}(C) = \frac{1}{2} 
\{(x_1,\ldots,x_n) \in \ZZ^n \mid
(x_1 \bmod 4,\ldots,x_n \bmod 4)\in C\}
\]
is an even (resp.\ odd) unimodular lattice with 
minimum norm $\min\{4,d_E/4\}$ \cite{Z4-BSBM}.
Let $C$ be a Type~I code of length $n$ and minimum Euclidean
weight $d_E$.
Since $A_4(C)$ is an odd unimodular lattice with
minimum norm $\min\{4,d_E/4\}$,
$\mu^O_{max}(n)$ in (\ref{Eq:mu})
implies upper bounds,
which improve the upper bounds described in 
Section \ref{Sec:1}, for $25 \le n \le 35, n\ne 32$.
Hence, we have the following bounds
\begin{equation}\label{Eq:dE}
d_{max,E}^I(n) \le
\left\{\begin{array}{ll}
8  & (n=25),                    \\
12 & (n=26,27,\ldots,31,33,34,35), \\
16 & (n=32, 36,37, \ldots, 47).
\end{array}\right.
\end{equation}

\section{Optimal Type~I $\ZZ_4$-codes}

In this section, 
we determine the largest minimum Euclidean weights 
$d_{max,E}^I(n)$.
This is done by constructing Type~I codes
meeting the bound (\ref{Eq:dE}) with equality
for lengths $n=25,\ldots,29,31,\ldots,36$, $41,43,44,45$.

\begin{itemize}
\item Lengths $26, \ldots, 29, 33, 34, 36, 41, 43,44,45$:

We have found a binary doubly even code $B_n$ of length $n$ with
$d(B_n)=12$ and $d(B_n^\perp)\ge3$ for $n=26, \ldots, 29, 33, 34$, and
$d(B_n)=16$ and $d(B_n^\perp)\ge4$ for $n=36,41,43,44,45$.
Note that there is a self-dual $\ZZ_4$-code $C$ with $C^{(1)}=B$
for any given binary doubly even code $B$ (see \cite{Z4-PLF}).
It follows from (\ref{eq:dE}) that there is a Type~I code 
meeting the bound (\ref{Eq:dE}) with equality for these lengths.

The method of construction of a self-dual $\ZZ_4$-code $C$ 
with $C^{(1)}=B$ was given in~\cite[Section 3]{Z4-PLF}.
Using this method, we explicitly have found an optimal Type~I code $C_n$
with $C_n^{(1)}=B_n$
for these lengths.
In order to save space, 
instead of listing generator matrices, 
we only list in Figure \ref{Fig:1}
the $k_1 \times (n-k_1)$ matrix
\[
M_n=
\left(\begin{array}{cc}
A & B_1+2B_2 \\
\end{array}\right),
\]
in standard form (\ref{Z4-g-matrix}), 
since the lower part of (\ref{Z4-g-matrix}) can be
obtained from $M_n$ for each code $C_n$.
The minimum Euclidean, Lee, Hamming weights
$d_E,d_L,d_H$ of $C_n$ 
(see e.g.\ \cite{Z4-BSBM} for the definition of $d_L$) 
are listed in  Table \ref{Table}.
The parameters $[n,k,d]$ and the orders $\#\Aut$ of the automorphism groups
of their residue codes $C_n^{(1)}$
are also listed in Table \ref{Table}.
Note that the minimum Hamming weight of a self-dual
$\ZZ_4$-code $C$ is the same as
$d(C^{(2)})$ \cite{Rains}.

%

\begin{figure}[p]
\centering
{\footnotesize
\begin{tabular}{ll}
$
M_{26}=
\left(\begin{array}{c}
01111011110001230212\\
11001011111000101120\\
11011110000100311232\\
00100110101110333302\\
11000110111101210021\\
01110001111110003003
\end{array}\right)
$
&
$
M_{27}=
\left(\begin{array}{c}
10011100010113320130\\
01011010100013333030\\
01010101010002131131\\
00010010111112303103\\
11110010010101330201\\
11111111110002220230\\
01111100101111000003
\end{array}\right)
$
\\ 
$
M_{28}=
\left(\begin{array}{c}
001101001010102330111\\
000110100101011211211\\
100011010010101123101\\
010001101001011132130\\
101000110100100331031\\
110100001010011211101\\
011010010101001121332
\end{array}\right)
$
&
$
M_{29}=
\left(\begin{array}{c}
1110001000001111103030\\
1101100110101002001303\\
1100001001101112030103\\
1101110010010002110213\\
0110111000111010010212\\
1010100110011011202011\\
0001110111111003232200
\end{array}\right)
$
\\ 
$
M_{33}=
\left(\begin{array}{c}
100101011010001020213213\\
110010100011100203101010\\
011111010001001010333131\\
110111000110111200023000\\
001010111111011113223001\\
101011111111100312212322\\
001001111100011302030032\\
111110101111110230220033\\
010110111101001020201010
\end{array}\right)
$
&
$
M_{34}=
\left(\begin{array}{c}
010110100101010311033313\\
001001110101000113202213\\
011010010101003100021323\\
000111011101110221331112\\
101000101001001031101030\\
011111010001100322300012\\
110010100111111103122323\\
111010010111112213011100\\
110000111011110000232003\\
111111111100002202022023
\end{array}\right)
$
\\
$
M_{36}=
\left(\begin{array}{c}
11100000110011011001003323121\\
01110110011000100010011231213\\
11010011111010100111111102121\\
01000011111101001010110130212\\
00011010110110010100013213121\\
00011100011010111101111022010\\
11111111111111000000000002023
\end{array}\right)
$
&
$
M_{41}=
\left(\begin{array}{c}
1001101011010100100110113101331\\
1011111100100111100012210033331\\
1100110001000111100003301002331\\
0001000111100111010012301003121\\
0011111001111011101002203222032\\
1000011011111010001111023333310\\
1011111111100101011011201310220\\
0110100100000101010013122113213\\
1000010101101100100100130112332\\
1110110011011110111101132113031
\end{array}\right)
$
\end{tabular}
\caption{Generator matrices}
\label{Fig:1}
}
\end{figure}

\setcounter{figure}{0}
\begin{figure}[htb]
\centering
{\scriptsize
\begin{tabular}{ll}
$
M_{43}= 
\left(\begin{array}{c}
101010101001101010001001103210230\\
110011001110110011000012211331132\\
001001011111011101100013131220013\\
101101111010010011001003033110211\\
010101110000001100110003300330330\\
111111101011001001111101102202120\\
101101111001101111001112230012100\\
010011110011110011011111020211320\\
000100111000111110111113022221133\\
111111111111110000000000020232222
\end{array}\right)
$
 &
$
M_{44}=
\left(\begin{array}{c}
011110001001001001011010122113202\\
001111000100100100101103210031100\\
000111100011010010010120101021330\\
100011110001101001001002010320313\\
010001111000110100100130203010013\\
001000111101011010010011022123003\\
000100011110101101001031320012300\\
100010001110010110100103312003012\\
110001000111001011010002331022301\\
111000100010100101101032013122012\\
111100010000010010110121023132003
\end{array}\right)
$
\\ 
$
M_{45}=
\left(\begin{array}{c}
000110111110100101100000111131021032\\
100011011111010010110000011213322301\\
110001101011101001111000001121332032\\
111000110101110100111100000232133223\\
011100011010111010011110000323013320\\
101110001001011101001111000232103112\\
110111000100101110000111100001230333\\
011011100010010111000011110120101213\\
001101110101001011000001111130230121
\end{array}\right)
$
\end{tabular}
\caption{Generator matrices (continued)}
\label{Fig:2}
}
\end{figure}

\begin{table}[thb]
\caption{Optimal Type~I $\ZZ_4$-codes}
\label{Table}
\begin{center}
{\small
\begin{tabular}{c||c|c|c||c|c}
\noalign{\hrule height0.8pt}
& \multicolumn{3}{c||}{$C_n$} & \multicolumn{2}{c}{$C_n^{(1)}$} \\
\hline
Code & $d_E$ & $d_L$ & $d_H$ & $[n,k,d]$ & $\#\Aut$ \\
\hline
$C_{26}$ & 12 & 6 & 3 & $[26,  6, 12 ]$& 120 \\
$C_{27}$ & 12 & 6 & 3 & $[27,  7, 12 ]$& 240 \\
$C_{28}$ & 12 & 8 & 4 & $[28,  7, 12 ]$& 10752 \\
$C_{29}$ & 12 & 6 & 3 & $[29,  7, 12 ]$& 1 \\
$C_{33}$ & 12 & 6 & 3 & $[33,  9, 12 ]$& 1 \\
$C_{34}$ & 12 & 6 & 3 & $[34, 10, 12 ]$& 1 \\
$C_{36}$ & 16 & 8 & 4 & $[36,  7, 16 ]$& 1451520 \\
$C_{41}$ & 16 & 8 & 4 & $[41, 10, 16 ]$& 1 \\
$C_{43}$ & 16 & 8 & 4 & $[43, 10, 16 ]$& 1 \\
$C_{44}$ & 16 & 8 & 4 & $[44, 11, 16 ]$& 11 \\
$C_{45}$ & 16 & 8 & 4 & $[45,  9, 16 ]$& 9 \\
\noalign{\hrule height0.8pt}
   \end{tabular}
}
\end{center}
\end{table}


\item {Length $32$}:

Let $L$ be a unimodular lattice in dimension $n$
and let $k$ be a positive integer.
A set $\{f_1, \ldots, f_{n}\}$ of $n$ vectors 
$f_1, \ldots, f_{n}$ in $L$ with
$(f_i, f_j) = k \delta_{ij}$
is called a {\em $k$-frame} of $L$,
where $\delta_{ij}$ is the Kronecker delta.
It is known that an even (resp.\ odd) unimodular lattice $L$ 
contains a 4-frame if and 
only if there is a Type~II (resp.\ Type~I)
$\ZZ_4$-code $C$ with $A_4(C) \cong L$.

Conway and Sloane \cite{CS98} showed that there are
exactly five odd unimodular lattices in dimension $32$
having minimum norm $4$, up to isomorphism.
In addition, 
such a lattice $\Lambda$
contains vectors of the form
\[
\frac{1}{\sqrt{8}}(\pm 4, \pm 4, 0, \ldots,0),
\ldots,
\frac{1}{\sqrt{8}}(0, \ldots,0, \pm 4, \pm 4).
\]
Hence, $\Lambda$ contains a $4$-frame.
This means that there is a Type~I $\ZZ_4$-code $C_{32}$
of length $32$ with $A_4(C_{32}) \cong \Lambda$.
Since $\Lambda$ has minimum norm $4$,
$C_{32}$ has minimum Euclidean weight $16$, that is,
$C_{32}$ is optimal.

\item{Lengths $25, 31, 35$}:

Let $C$ be a self-dual code of length $n$ $(n \ge 2)$.
Then the following code
\[
\sub(C) = \{(x_2,\ldots,x_{n}) \mid
(x_1,x_2,\ldots,x_n) \in C, x_1 \in \{0,2\}\}
\]
is a self-dual code of length $n-1$.
The codes $\sub(C_{26})$, $\sub(C_{32})$ and $\sub(C_{36})$
are self-dual codes of
lengths $25,31$ and $35$, respectively.
Moreover, from (\ref{Eq:dE}),
the codes $\sub(C_{26})$, $\sub(C_{32})$ and $\sub(C_{36})$
have minimum Euclidean weights $8,12$ and $12$, respectively,
since $C_{26}$, $C_{32}$ and $C_{36}$
have minimum Euclidean weight
$12,16$ and $16$, respectively.

Since there are self-dual codes of lengths $8, 17$
and minimum Euclidean weight $8$, 
the direct sum of the codes is also a self-dual code
of length $25$ and minimum Euclidean weight $8$.

\end{itemize}

By constructing Type~I codes 
meeting the bound (\ref{Eq:dE}) with equality
for lengths $n=25,\ldots,29,31,\ldots,36$, $41,43,44,45$,
we determine the largest minimum Euclidean weight 
$d_{max,E}^I(n)$ for $n \le 47,\ n\not=37$,
as follows.

\begin{prop}
Let $d_{max,E}^I(n)$ denote
the largest minimum Euclidean weight among 
Type~I $\ZZ_4$-codes of length $n$.
Then
$d_{max,E}^I(25)=8$,
$d_{max,E}^I(n)=12$ if $n=26,\ldots,31,33,34,35$, and 
$d_{max,E}^I(n)=16$ if $n=32,36,38,\ldots,47$.
\end{prop}

For length $37$,
our extensive search failed to discover 
a Type~I $\ZZ_4$-code with minimum Euclidean weight $16$.
However, 
we have found a Type~I code 
with minimum Euclidean weight $12$.
Hence, $d_{max,E}^I(37)=12$ or $16$
(see Table \ref{Tab:dE}).

\section{Optimal odd unimodular lattices}

\subsection{Dimension 41}

By Construction A, optimal Type~I $\ZZ_4$-codes $C_n$
constructed in the previous section give
optimal odd unimodular lattices $A_4(C_{n})$.
In particular, 
the first explicit example of an optimal
odd unimodular lattice $A_4(C_{41})$ in dimension $41$ 
can be constructed from $C_{41}$.

\begin{prop}
There is an odd unimodular lattice in dimension $41$
having minimum norm $4$.
\end{prop}


We consider the theta series of optimal odd unimodular lattices
in dimension $41$.
Conway and Sloane \cite{CS98} show that
if the theta series of an odd unimodular lattice $L$
in dimension $n$
is written as 
\begin{equation}
\label{Eq:theta}
\theta_L(q)=
 \sum_{j =0}^{\lfloor n/8\rfloor} a_j\theta_3(q)^{n-8j}\Delta_8(q)^j,
\end{equation}
then the theta series of the shadow $S$ 
(see \cite{CS98} for the definition)
is written as
\begin{equation}
\label{Eq:theta-S}
\theta_S(q)= \sum_{j=0}^{\lfloor n/8\rfloor}
\frac{(-1)^j}{16^j} a_j\theta_2(q)^{n-8j}\theta_4(q^2)^{8j}
= \sum_i B_i q^i \text{ (say)},
\end{equation}
where 
$\Delta_8(q) = q \prod_{m=1}^{\infty} (1 - q^{2m-1})^8(1-q^{4m})^8$
and $\theta_2(q), \theta_3(q)$ and $\theta_4(q)$ are the Jacobi 
theta series \cite{SPLAG}.
As the additional conditions, it holds that
there is at most one nonzero $B_r$ for $r < (\mu+2)/2$;
$B_r=0$ for $r < \mu/4$; and $B_r \le 2$ for $r < \mu/2$,
where $\mu$ is the minimum norm of $L$.

In the case $n=41$, since minimum norm $\mu$ is $4$,
$a_0, \ldots, a_3$ in (\ref{Eq:theta}) and (\ref{Eq:theta-S})
are determined as follows:
$a_0=1$, $a_1=-82$, $a_2=1476$ and $a_3=-3280$.
Since the coefficients in the shadow
must be non-negative integers,
$a_4$ is divisible by $2^7$ and 
$a_5$ is divisible by $2^{19}$.
Thus, we put $a_4=2^7 \alpha$ and $a_5=-2^{19} \beta$.
Then we have the possible theta series $\theta_{L}$ and
$\theta_S$ of an optimal odd unimodular lattice $L$ 
in dimension $41$ and its shadow $S$:
\begin{align*}
\theta_L &=
1 
+ (15170 + 128\alpha )q^4 
+ (1226720 - 1792\alpha - 524288\beta) q^5 
\\ & \hspace*{3cm}
+ (42928640 + 8192\alpha + 19922944\beta) q^6 
+ \cdots \text{ and}
\\
\theta_S &=
\beta q^{1/4} 
+ (\alpha - 79\beta) q^{9/4} 
+ (104960 - 55\alpha + 3040\beta)q^{17/4} 
+ \cdots,
\end{align*}
respectively, 
where $\beta=0$ or $\alpha=79\beta$ by the above additional conditions.
By calculating the kissing number of $A_4(C_{41})$ and 
the minimum norm of its shadow,
we determine the theta series of the lattice $A_4(C_{41})$ as
follows:
\begin{multline*}
1 + 15426 q^4 + 1223136 q^5 + 42945024 q^6 + 867179520 q^7 
\\ 
+ 11719744560 q^8 + 116521216256 q^9 + 909236984832 q^{10}
+ \cdots.
\end{multline*}

\subsection{Minimum norms and kissing numbers}
In Table \ref{Tab:L},
we list the minimum norms $\mu(L)$ 
 and the kissing numbers $N(L)$
of optimal odd unimodular lattices $L=A_4(C_{n})$ 
constructed from $C_n$ given in Table \ref{Table}.

\begin{table}[thb]
\caption{Minimum norms and kissing numbers}
\label{Tab:L}
\begin{center}
{\small
\begin{tabular}{c|cc||c|cc}
\noalign{\hrule height0.8pt}
$L$ &  $\mu(L)$ & $N(L)$ &
$L$ &  $\mu(L)$ & $N(L)$ \\
\hline
$A_4(C_{26})$ & 3 &  3120 & $A_4(C_{36})$ & 4 & 51032 \\
$A_4(C_{27})$ & 3 &  2664 & $A_4(C_{41})$ & 4 & 15426 \\
$A_4(C_{28})$ & 3 &  1728 & $A_4(C_{43})$ & 4 &  9286 \\
$A_4(C_{29})$ & 3 &  1856 & $A_4(C_{44})$ & 4 &  8392 \\
$A_4(C_{33})$ & 3 &   752 & $A_4(C_{45})$ & 4 &  7866 \\
$A_4(C_{34})$ & 3 &   528 &          && \\
\noalign{\hrule height0.8pt}
   \end{tabular}
}
\end{center}
\end{table}

For dimensions up to $28$, optimal odd unimodular lattices
have been classified (see \cite[p.~xliii--xliv]{SPLAG}).
Borcherds \cite{Bor} showed that
there is a unique optimal odd unimodular lattice $S_{26}$
in dimension $26$ (see \cite[p.~xliii]{SPLAG}).
Hence, the lattice $A_4(C_{26})$ gives 
an alternative construction of $S_{26}$.
We list the symmetrized weight enumerator $swe_{26}$ (see \cite{Z4-CS}
for the definition) of $C_{26}$
at the end of this section.
Bacher and Venkov \cite{BV} showed that
there are three (resp.\ $38$)
non-isomorphic optimal odd unimodular
lattices in dimension $27$ (resp.\ $28$)
(see \cite[p.~xliv]{SPLAG}).
By comparing the kissing numbers and the automorphism groups,
we have that $A_4(C_{27}) \cong 
{\mathbf R}_{\bf 27,1}(\emptyset)$ in \cite[Table 4]{BV},
and $A_4(C_{28}) \cong
{\mathbf R}_{\bf 28,38e}(\emptyset)$ in \cite[Table 5]{BV}.
Hence, these lattices $S_{26}$, 
${\mathbf R}_{\bf 27,1}(\emptyset)$ and 
${\mathbf R}_{\bf 28,38e}(\emptyset)$
contain a $4$-frame.
We remark that the lattices have no $3$-frame.

For other dimensions, 
since the lattices $A_4(C_n)$ ($n=33,36,44$) have
different kissing numbers than those of 
the known lattices in \cite{S-http}
and no example was given in \cite{S-http} for dimension $30$,
the lattices $A_4(C_n)$ ($n=30, 33,36,44$) 
provide other examples of optimal odd unimodular lattices.

Since an odd unimodular lattice in dimension $41$
having minimum norm $4$ has been constructed,
the largest minimum norm $\mu^O_{max}(41)$ is $4$.
The smallest dimension $n$ for which the largest minimum
norm $\mu^O_{max}(n)$ has not been determined is $37$.
Hence, it is worthwhile to determine if there is a Type~I $\ZZ_4$-code of 
length $37$ and minimum Euclidean weight $16$.
At dimension $48$, the largest minimum norm $\mu^O_{max}(48)$ is 
exactly $5$ 
(\cite{Gaulter}, \cite{Venkov} and \cite{Rains-Sloane98}).
No Type~I $\ZZ_4$-code constructs an odd unimodular 
lattice with minimum norm $5$ by Construction A.
It seems that the connection between self-dual $\ZZ_4$-codes and
unimodular lattices is no longer useful at this point.

{\footnotesize
\begin{align*}
swe_{26}=
&
a^{26} + 30 a^{23} c^3 + 255 a^{22} c^4 + 1100 a^{21} c^5 
+ 3571a^{20} c^6 + 9990 a^{19} c^7 + 24330 a^{18} c^8 
\\ &
+ 49680 a^{17} c^9 + 83237 a^{16} c^{10} + 119004 a^{15} c^{11} 
+ 2880 a^{14} b^{12} + 150750 a^{14} c^{12} 
\\ &
+ 40320 a^{13} b^{12} c + 164680 a^{13} c^{13} + 262080 a^{12} b^{12} c^2 
+ 150750 a^{12} c^{14} + 1048320 a^{11} b^{12} c^3 
\\ &
+ 119004 a^{11} c^{15} + 17408 a^{10} b^{16} + 2882880 a^{10} b^{12} c^4 
+ 83237 a^{10} c^{16} + 174080 a^9 b^{16} c 
\\ &
+ 5765760 a^9 b^{12} c^5 
+ 49680 a^9 c^{17} + 783360 a^8 b^{16} c^2 + 8648640 a^8 b^{12} c^6 
+ 24330 a^8 c^{18} 
\\ &
+ 2088960 a^7 b^{16} c^3 + 9884160 a^7 b^{12} c^7 
+ 9990 a^7 c^{19} + 16384 a^6 b^{20} + 3655680 a^6 b^{16} c^4
\\ &
+ 8648640 a^6 b^{12} c^8 + 3571 a^6 c^{20} + 98304 a^5 b^{20} c 
+ 4386816 a^5 b^{16} c^5 + 5765760 a^5 b^{12} c^9 
\\ &
+ 1100 a^5 c^{21} 
+ 245760 a^4 b^{20} c^2 + 3655680 a^4 b^{16} c^6 
+ 2882880 a^4 b^{12} c^{10} + 255 a^4 c^{22} 
\\ &
+ 327680 a^3 b^{20} c^3 + 2088960 a^3 b^{16} c^7 
+ 1048320 a^3 b^{12} c^{11} + 30 a^3 c^{23} + 245760 a^2 b^{20} c^4 
\\ &
+ 783360 a^2 b^{16} c^8 + 262080 a^2 b^{12} c^{12} + 98304 a b^{20} c^5 
+ 174080 a b^{16} c^9  + 40320 a b^{12} c^{13} 
\\ &
+ 16384 b^{20} c^6 + 17408 b^{16} c^{10} 
+ 2880 b^{12} c^{14} + c^{26}
\end{align*}
}

\bigskip
\noindent {\bf Acknowledgments.}
The author would like to thank the anonymous referees for 
their helpful comments.


\end{document}